\documentclass[11pt,a4paper]{article}
\usepackage[active]{srcltx}
\usepackage{amsfonts}
\usepackage{amsmath}
\usepackage{graphicx,epsfig}
\usepackage{amssymb}
\usepackage{theorem,subfigure}
\usepackage{rotating}
\usepackage[usenames,dvipsnames]{color} 
\usepackage[small,bf]{caption}
\setcaptionwidth{14.0cm}

\newtheorem{theorem}{Theorem}[section]

\newtheorem{algorithm}[theorem]{Algorithm}

\newtheorem{corollary}[theorem]{Corollary}

\newcommand{\IN}{\mbox{$\mathbb N$}}
\newcommand{\IR}{\mbox{$\mathbb R$}}

\def\varddots{\mathinner{\mkern1mu\raise1pt\vbox{\kern7pt\hbox{.}}\mkern2mu
    \raise4pt\hbox{.}\mkern2mu\raise7pt\hbox{.}\mkern1mu}}
\font\mfett=cmmib10 at10pt
\def\nufett{\hbox{\mfett\char023}}

\usepackage{marginnote}
\usepackage{ulem}
\usepackage[top=5cm,left=2.5cm, heightrounded,marginparsep=-0.5cm, marginparwidth=4.5cm]{geometry}

\begin{document}
\title {Pseudo-inverses of difference matrices and their application to sparse signal approximation}
\author{Gerlind Plonka\footnote{Institute for Numerical and Applied Mathematics,
University of G\"ottingen,
Lotzestr. 16--18, 37083 G\"ottingen, Germany.Email: plonka@math.uni-goettingen.de}
\qquad Sebastian Hoffmann\footnote{Mathematical Image Analysis Group,
Faculty of Mathematics and Computer Science, Campus E1.7,
Saarland University, 66041 Saarbr\"ucken, Germany. Email: hoffmann@mia.uni-saarland.de} \qquad
Joachim Weickert\footnote{Mathematical Image Analysis Group,
Faculty of Mathematics and Computer Science, Campus E1.7,
Saarland University, 66041 Saarbr\"ucken, Germany. Email: weickert@mia.uni-saarland.de}}

\maketitle
\abstract{
We derive new explicit expressions for the components of Moore-Penrose inverses of symmetric difference matrices.
These generalized inverses are applied in a new regularization approach for scattered data interpolation based on partial differential equations. The columns of the Moore-Penrose inverse then serve as elements of a dictionary that allow a sparse signal approximation. In order to find a set of suitable data points for signal representation we apply the orthogonal patching pursuit (OMP) method. 
}

\section{Introduction}

Within the last years, there have been several attempts to derive new sparse representations for signals and images in the context of sparse data interpolation with partial differential equations.
The essential key of these so called inpainting methods is to fix a suitable set of data points  that admits  a very good signal or image approximation, when
the unspecified data are interpolated by means of a diffusion process.
Finding optimal sparse data leads to a challenging optimization problem that is NP-hard, and different heuristic strategies have been proposed.
Using  nonlinear anisotropic diffusion processes, this  approach has been introduced in \cite{G05} and strongly improved since that time by more sophisticated  optimization methods such as \cite{SPM14}.
Inpainting  based on linear differential operators such as the Laplacian is 
conceptually simpler, but can still provide sparse signal approximations 
with high quality; see e.g.\ \cite{HSW13,MHW12,OCBP14}. These recent results 
outperform earlier attempts in this direction that use specific features such 
as edges \cite{Ca88,HM89,El99} or toppoints in scale-space \cite{JSGA86}. 
\medskip

In this paper, we focus on the one-dimensional discrete case. The goal is to provide new insights into this problem from a linear algebra perspective.
In the first part  we derive a formula for the Moore-Penrose inverse of difference matrices. This formula has already been proposed for finding generalized inverses of circulant matrices of size $N \times N$ with rank $N-1$ in \cite{S79} and for discrete Laplace operators considered e.g.\ in \cite{XG03,RZB13}.
We generalize this formula to arbitrary symmetric $N \times N$ difference matrices of rank $N-1$ possessing the eigenvector ${\bf 1} = (1,\ldots , 1)^{T}$ to the single eigenvalue $0$. This formula enables us to derive explicit expressions for the  components of the Moore-Penrose inverse for special difference matrices that play a key role in discrete diffusion inpainting. In particular, we study difference matrices that approximate the second and fourth order signal derivative, i.e., the one-dimensional Laplace operator
and the biharmonic operator, with periodic resp.\ reflecting (homogeneous Neumann) boundary conditions.

With the help of the generalized inverse of the difference matrices, we propose a new regularization approach that admits different optimization strategies for finding a suitable sparse set of data points for signal reconstruction. In particular, we employ the orthogonal matching pursuit (OMP) method as a conceptually simple and efficient greedy algorithm to construct the desired set of data points. 
\medskip

The paper is organized as follows.
In Section 2 we present the formula for the Moore-Penrose inverse of symmetric difference matrices and derive
  explicit  expressions for the components of the Moore-Penrose inverse of difference matrices 
corresponding to the discrete Laplace and biharmonic operator using the theory of difference equations.
Section 3 is devoted to a new regularization approach to the discrete inpainting problem.
The columns of the explicitly given Moore-Penrose inverses of the difference matrices can be understood as discrete Green's functions \cite{HPW15} and  are used as a dictionary for the orthogonal matching pursuit algorithm (OMP). Since the original discrete inpainting problem (with the Laplace or biharmonic operator) is equivalent to  a least squares spline approximation problem with free knots, 
our approach  can also be seen  as an alternative method to solve this nonlinear approximation problem. 
Finally we present some numerical results in Section 4.

\section{Explicit Moore-Penrose inverses for difference matrices}
\label{sec2}

Let ${\bf L}$ be an $N \times N$ matrix. We say  that 
${\bf L}$ is a {\it difference matrix} if it satisfies the following properties.
\begin{description}
\item{(i)} The matrix ${\bf L}$  is symmetric and of rank $N-1$;
\item{(ii)}  ${\bf L} {\bf 1} = {\bf 0}$ where ${\bf 1} := (1, \ldots , 1)^{T} \in {\IR}^{N}$, i.e., ${\bf 1}$ is the  eigenvector of ${\bf L}$ to the single eigenvalue $0$.
\end{description}
 
Matrices of this type typically occur when discrete approximations of symmetric differential operators with periodic or Neumann boundary conditions are considered. 

In particular, a discretization of the second derivative (one-dimensional Laplace) with periodic boundary conditions yields a circulant matrix ${\bf L}$ of the form 
\begin{equation}\label{A1} {\bf L} = {\bf A} = {\rm circ} \,  (-2,1,0, \ldots , 0, 1)^{T} :=\left( \begin{array}{ccccc} 
-2 & 1 & 0 &  \ldots & 1 \\
1 & -2 & 1 &  &   \\
0 & \ddots & \ddots & \ddots & \\
 &  & \ddots & \ddots & 1 \\
 1 & 0 & \ldots  & 1  & -2  \end{array} \right) \in {\IR}^{N \times N},
\end{equation}
while for Neumann conditions, we consider
\begin{equation}\label{A2} {\bf L} = {\bf B}  :=\left( \begin{array}{ccccc} 
-1 & 1 & 0 &  \ldots & 0 \\
1 & -2 & 1 &  &   \\
0 & \ddots & \ddots & \ddots & \\
 &  & \ddots & -2 & 1 \\
 0 & 0 & \ldots  & 1  & -1  \end{array} \right) \in {\IR}^{N \times N}.
\end{equation}
As a  discretization of the fourth derivative (one-dimensional biharmonic operator) we obtain in the periodic case
$$ {\bf L} = -{\bf A}^{2} = {\rm circ} \,  (-6,4,-1,0 \ldots , 0, -1,4)^{T} $$
and in the case of reflecting boundary conditions
$$ {\bf L} = - {\bf B}^{2}. $$
All these matrices satisfy the properties (i) and (ii) given above.

The following theorem generalizes the ideas of \cite{S79,B81} for pseudo-inverses of circulant matrices of rank $N-1$, as well as of \cite{XG03,RZB13} for pseudo-inverses of the Laplacian.
According to the usual definition, we say that ${\bf L}^{+}$ is the Moore-Penrose inverse of ${\bf L}$ if 
$$ {\bf L} {\bf L}^{+} {\bf L} = {\bf L}, \qquad  {\bf L}^{+} {\bf L} {\bf L}^{+} = {\bf L}^{+}, $$
and ${\bf L} {\bf L}^{+}$ as well as ${\bf L}^{+} {\bf L}$ are symmetric. With these properties, the Moore-Penrose inverse is uniquely defined; see \cite{P55}.

 \begin{theorem} \label{theoA+}
 For a difference matrix  ${\bf L} \in {\IR}^{N \times N}$ satisfying the conditions {\rm (i)} and {\rm (ii)} we obtain the Moore-Penrose inverse in the form
 \begin{equation} \label{A+} {\bf L}^{+} = ({\bf L} - \tau {\bf J})^{-1} + \frac{1}{\tau N^{2}} {\bf J},
 \end{equation}
 where $\tau \in {\IR} \setminus \{ 0 \} $ can be chosen arbitrarily, and ${\bf J}$ is the $N \times N$ matrix where every element is set to $1$.
 In particular, ${\bf L}^{+}$ does not depend on the choice of $\tau$. Further, ${\bf L}^{+}$ is symmetric.
 \end{theorem}
 
{\bf Proof}. First, we show that ${\bf L} - \tau {\bf J}$ is a regular matrix for all $\tau \in {\IR} \setminus \{ 0 \}$. Since ${\bf L}$ is symmetric, there  exists an orthonormal basis  of eigenvectors ${\bf v}_{0}, \ldots , {\bf v}_{N-1}$ for ${\bf L}$, where ${\bf v}_{0} = \frac{1}{\sqrt{N}}{\bf 1}$ is the eigenvector of ${\bf L}$ to the single eigenvalue $\lambda_{0}=0$. Hence, 
${\bf L} {\bf v}_{j} = \lambda_{j} {\bf v}_{j}$, $j=1, \ldots , N-1$, with eigenvalues $\lambda_{j} \in {\IR} \setminus \{ 0 \}$ implies 
$$ ({\bf L} - \tau {\bf J}) {\bf v}_{j} = \lambda_{j} {\bf v}_{j}$$
since ${\bf 1}^{T} {\bf v}_{j} =0$ for all $j=1, \ldots , N-1$.
Further, 
\begin{equation} \label{one}
 ({\bf L} - \tau {\bf J}) {\bf 1} = -\tau {\bf J} {\bf 1} = -\tau  N {\bf 1},
 \end{equation}
  i.e., $({\bf L} - \tau {\bf J})$ possesses the $N-1$ nonzero eigenvalues  $\lambda_{1}, \ldots , \lambda_{N-1}$ of ${\bf L}$ and the additional eigenvalue $-\tau N$. Therefore, it is regular. In order to show that ${\bf L}^{+}$ in (\ref{A+}) is the Moore-Penrose inverse of ${\bf L}$, we check the four properties that  determine tie Moore-Penrose inverse uniquely. From (\ref{one})
it follows that 
\begin{equation} \label{two}  ( {\bf L} - \tau {\bf J})^{-1} {\bf 1} = - \frac{1}{\tau  N} {\bf 1} \end{equation}
and hence
$$ ( {\bf L} - \tau {\bf J})^{-1} {\bf J} = - \frac{1}{\tau  N} {\bf J}. $$
Analogously, we simply derive ${\bf J} ({\bf L}- \tau {\bf J})^{-1} = -\frac{1}{\tau N
} {\bf J}$ from ${\bf 1}^{T} ({\bf L} - \tau {\bf J}) = - \tau N {\bf 1}^{T}$. Using ${\bf J} {\bf L} = {\bf L} {\bf J} = {\bf 0}$ we find 
\begin{eqnarray}
\nonumber
{\bf L} {\bf L}^{+} {\bf L} &=& {\bf L} \left( ({\bf L} - \tau {\bf J
})^{-1} + \frac{1}{\tau N^{2}} {\bf J} \right) {\bf L} \\
\nonumber
&=& {\bf L} ({\bf L} - \tau {\bf J})^{-1} ({\bf L} - \tau {\bf J} + \tau {\bf J}) \\
\nonumber
&=& {\bf L} ({\bf I} + ({\bf L} -\tau {\bf J})^{-1} \tau {\bf J})\\
&=& {\bf L} \left( {\bf I} - \frac{\tau}{\tau N} {\bf J}\right)  
\label{equationref1} \\ 
\nonumber
&=& {\bf L}, 
\end{eqnarray}
and by ${\bf J}^{2} = N {\bf J}$, 
\begin{eqnarray}
\nonumber
{\bf L}^{+} {\bf L} {\bf L}^{+} 
&=& {\bf L}^{+} {\bf L} \left( ({\bf L} - \tau {\bf J})^{-1} +
\frac{1}{\tau N^{2}} {\bf J}\right) \\
\nonumber
&=& {\bf L}^{+} ({\bf L}- \tau {\bf J} + \tau {\bf J}) 
\left( ({\bf L} - \tau {\bf J})^{-1} + \frac{1}{\tau N^{2}} {\bf J}\right) \\
\nonumber
&=& {\bf L}^{+} \left( {\bf I} + \frac{1}{\tau N^{2}} ({\bf L} - \tau {\bf J}) {\bf J} + \tau {\bf J} ({\bf L} - \tau {\bf J})^{-1} + \frac{\tau}{\tau  N^{2}} {\bf J}^{2} \right) \\
\nonumber
&=& {\bf L}^{+} \left({\bf I} - \frac{\tau}{\tau N^{2}} {\bf J}^{2} - \frac{\tau}{\tau N} {\bf J} + \frac{1}{N} {\bf J} \right)\\
&=& {\bf L}^{+} \left({\bf I} - \frac{1}{ N} {\bf J}\right)
\label{equationref2} \\ 
\nonumber
&=& {\bf L}^{+} - \frac{1}{ N} \left(({\bf L} - \tau {\bf J})^{-1} + \frac{1}{\tau N^{2}} {\bf J}\right) {\bf J} \\
\nonumber
&=& {\bf L}^{+} - \frac{1}{ N} \left( -\frac{1}{\tau N} + \frac{1}{\tau N}\right) {\bf J} \\
\nonumber
&=& {\bf L}^{+}.
\end{eqnarray}
The equalities (\ref{equationref1}) and (\ref{equationref2}) particularly imply ${\bf L} {\bf L}^{+} = {\bf L}^{+} {\bf L}= ({\bf I} - \frac{1}{N} {\bf J})$, i.e., ${\bf L}^{+} {\bf L}$ as well as ${\bf L} {\bf L}^{+}$ are symmetric. 
Since ${\bf L}^{+}$ is uniquely defined, it does not depend on the choice of $\tau$. Finally, we have
\[\left({\bf L}^{+}\right)^\top = \left({\bf L}^\top\right)^{+} = {\bf L}^{+},\]
which shows the symmetry of ${\bf L}^{+}$.\hfill $\Box$
\medskip

We want to apply (\ref{A+}) in order to derive explicit expressions for the special difference matrices ${\bf A}$ and ${\bf B}$ in (\ref{A1}) and (\ref{A2}) as well as for $-{\bf A}^{2}$ and $-{\bf B}^{2}$ that approximate  differential operators of second and fourth order, 
together with periodic or reflecting (homogeneous Neumann) boundary conditions.

\begin{theorem}\label{theo2}
The Moore-Penrose inverse of the circulant matrix ${\bf A}$ in {\rm (\ref{A1})} is a symmetric matrix of the form
$$ {\bf A}^{+} = {\rm circ} \, (a_{0}^{+}, a_{1}^{+}, \ldots , a_{N-1}^{+}) = \left( \begin{array}{ccccc}
a_{0}^{+} & a_{N-1}^{+} & a_{N-2}^{+} & \ldots & a_{1}^{+} \\
a_{1}^{+} & a_{0}^{+} & a_{N-1}^{+} & \ldots & a_{2}^{+} \\
a_{2}^{+} & a_{1}^{+} & a_{0}^{+} &  & a_{3}^{+} \\
\vdots &  & \ddots & \ddots & \vdots \\
a_{N-1}^{+} & a_{N-2}^{+} & \ldots & a_{1}^{+} & a_{0}^{+} \end{array}\right),
$$
where
$$ a_{j}^{+} = \frac{1}{12 N} (1-N^{2})  + \frac{1}{2N} j(N-j), \qquad j=0, \ldots, N-1.$$
\end{theorem}

{\bf Proof}.
By Theorem \ref{theoA+} with $\tau = -1$ we have ${\bf A}^{+} = ({\bf A} + {\bf J})^{-1} - \frac{1}{N^{2}} {\bf J}$, where ${\bf J}$ is the matrix with every element one. Since ${\bf A} + {\bf J}$ is circulant, also $({\bf A} + {\bf J})^{-1}$ is circulant, and we assume that it is of the form
$({\bf A} + {\bf J})^{-1} = {\rm circ} \, (b_{0}, b_{1}, \ldots , b_{N-1})$. The proof of Theorem \ref{theoA+} implies 
$${\bf A} ({\bf A} + {\bf J})^{-1} = ({\bf A} + {\bf J} - {\bf J}) ({\bf A} + {\bf J})^{-1} =  {\bf I} - {\bf J}({\bf A} + {\bf J})^{-1} = {\bf I} - \frac{1}{N} {\bf J}, $$ i.e., we obtain
\begin{equation}\label{b1}
b_{j-1}- 2b_{j} + b_{j+1} = - \frac{1}{N}, \qquad j=1, \ldots , N-1
\end{equation}
and the boundary conditions $b_{N-1} - 2 b_{0} + b_{1} = 1- \frac{1}{N}$
and $b_{N-2} - 2b_{N-1} + b_{0} = - \frac{1}{N}$. Further, by (\ref{two}) it follows that $({\bf A} + {\bf J})^{-1} {\bf 1} = \frac{1}{N} {\bf 1}$ and hence 
\begin{equation} \label{b2s} \sum_{j=0}^{N-1} b_{j} = \frac{1}{N}. 
\end{equation}
The inhomogenous difference equation (\ref{b1}) can be rewritten as
$$ \left( \begin{array}{c} b_{j} \\ b_{j+1} \end{array} \right) = 
\left( \begin{array}{cc} 0 & 1 \\ -1 & 2 \end{array} \right) \, 
\left( \begin{array}{c} b_{j-1} \\ b_{j} \end{array} \right) +
\left( \begin{array}{c} 0 \\ -\frac{1}{N} \end{array} \right), \qquad j=1, \ldots, N-2. $$
With ${\bf y}_{j+1} :=(b_{j}, b_{j+1})^{T}$, ${\bf g} := (0, -\frac{1}{N})^{T}$, and ${\bf G} := \left( \begin{array}{cc} 0 & 1 \\
-1 & 2 \end{array} \right)$ we obtain the linear, first order difference 
equation
$$ {\bf y}_{j+1} = {\bf G} {\bf y}_{j} + {\bf g} $$
with the general solution
$$ {\bf y}_{j} = {\bf G}^{j} {\bf y}_{0} + \sum_{\nu=1}^{j} {\bf G}^{j-\nu} {\bf g}, $$
where  ${\bf y}_{0}=(b_{0}, b_{1})^{T}$  needs to be determined by the boundary conditions and the additional condition (\ref{b2s}). A simple induction argument shows that the powers of ${\bf G}$ can be written as
$$ {\bf G}^{j} = \left( \begin{array}{cc} -j+1 & j \\ -j & j+1 \end{array} \right), $$
and we obtain
\begin{eqnarray} \nonumber
b_{j} &=& (1, 0) \, {\bf y}_{j} \\
\nonumber
&=&  (1,0 ) \left( \begin{array}{cc} -j+1 & j \\ -j & j+1 \end{array} \right)  \left( \begin{array}{c} b_{0} \\ b_{1} \end{array} \right) +
\sum_{\nu=1}^{j} (1,0 ) \left( \begin{array}{cc} -j+\nu +1 & j-\nu \\ \nu-j & j-\nu+1 \end{array} \right)  \left( \begin{array}{c} 0 \\ -\frac{1}{N} \end{array} \right) \\
\label{b2}
&=& b_{0} + j(b_{1} - b_{0}) - \frac{1}{N} \frac{j(j-1)}{2}.
\end{eqnarray}
Finally, we determine $b_{0}$ and $b_{1}$ from
\begin{equation}\label{b3}
 b_{N-1} - 2 b_{0} + b_{1} = 1 - \frac{1}{N}, \qquad \sum_{j=0}^{N-1} b_{j} = \frac{1}{N}
 \end{equation}
 while the second boundary condition is then automatically satisfied.
We insert (\ref{b2}) into the first equation of (\ref{b3}) and obtain
$$ b_{0} + (N-1)(b_{1} - b_{0}) - \frac{(N-2) (N-1)}{2 N} - 2 b_{0} + b_{1} = 1 - \frac{1}{N}, $$
i.e.,
$$ b_{1} - b_{0} = \frac{N-1}{2N}. $$
Hence, by (\ref{b2}), the second condition in (\ref{b3}) yields 
$$ \sum_{j=0}^{N-1} b_{j} = N b_{0} + \frac{N-1}{2N} \sum_{j=0}^{N-1} j
- \frac{1}{2N} \sum_{j=0}^{N-1} j(j-1) = \frac{1}{N}, $$
i.e.,
$$ N b_{0} = \frac{1}{N} - \frac{N^{2}-1}{12}. $$
With $b_{0} = \frac{1}{N^{2}} - \frac{N^{2}-1}{12N}$ we find
$$b_{j} = \frac{1}{N^{2}} + \frac{1-N^{2}}{12N} + \frac{j(N-j)}{2N}. $$
We finally obtain from (\ref{A+}) (with $\tau = -1$)
$$ a_{j}^{+} = b_{j} - \frac{1}{N^{2}} = \frac{1-N^{2}}{12N} + \frac{j(N-j)}{2N}.$$
Obviously ${\bf A}^{+}$ is symmetric and we have $a_{j}^{+} = a_{N-j}^{+}$.
\null \hfill $\Box$

Further, we can directly derive the Moore-Penrose inverse for the discrete biharmonic operator ${-\bf A}^{2}$ for periodic boundary conditions.

\begin{theorem}\label{theobi}
The Moore-Penrose inverse for the discrete biharmonic operator ${\bf L}=-{\bf A}^{2}$ with ${\bf A}$ in {\rm (\ref{A1})} with periodic boundary conditions is a symmetric matrix of the form
$$ {\bf L}^{+} = -({\bf A}^{2})^{+} = {\rm circ} \, (c_{0}^{+}, c_{1}^{+}, \ldots , c_{N-1}^{+}) $$
with $$ 
c_{k}^{+} = \frac{(1-N^{2} )(N^{2} +11)}{720 N} + \frac{k(N-k)(Nk-k^{2}+2)}{24N}, \qquad k=0, \ldots, N-1. $$
\end{theorem}

{\bf Proof}. From $({\bf A}^{2})^{+} = {(\bf A}^{+})^{2}$ it follows that
$$ - c_{k}^{+} =  \sum_{j=0}^{k} a_{j}^{+} a_{k-j}^{+} + \sum_{j=k+1}^{N-1} a_{j}^{+} a_{N+k-j}^{+}. $$
Inserting $a_{k}^{+} = \frac{1-N^{2}}{12N}
+ \frac{k(N-k)}{2N}$, a direct evaluation gives the desired result. 
Here one needs the Faulhaber formulas for sums of the
form $\sum_{j=1}^{r} j^{\ell} = \frac{B_{\ell+1}(r+1)-B_{\ell+1}(0)}{k+1}$, where $B_{\ell+1}$ denotes the Bernoulli polynomial of degree $\ell + 1$;
see e.g.\ \cite{K93}.
\hfill $\Box$
\medskip

As an interesting side information, we obtain the following trigonometric  identities.
\begin{corollary}
For $N \in {\IN}$ we have
\begin{eqnarray*}
\sum_{k=1}^{N-1} \frac{\cos(\frac{2\pi k j}{N})}{\sin^{2}(\frac{k\pi}{N})} &=& \frac{N^{2}-1}{3} - 2j(N-j), \qquad j=0, \ldots , N-1, \\
\sum_{k=1}^{N-1} \frac{\cos(\frac{2\pi k j}{N})}{\sin^{4}(\frac{k\pi}{N})} &=& \frac{(N^{2}-1)(N^{2}+11)}{45} - \frac{2j(N-j)(Nj-j^{2}+2)}{3}, \qquad j=0, \ldots , N-1, \\
\sum_{k=1}^{N-1} \frac{\sin(\frac{2\pi k j}{N})}{\sin^{2}(\frac{k\pi}{N})} &=& \sum_{k=1}^{N-1} \frac{\sin(\frac{2\pi k j}{N})}{\sin^{4}(\frac{k\pi}{N})} =0, 
\qquad j=0, \ldots , N-1,
\end{eqnarray*}
and particularly
$$ \sum_{k=1}^{N-1} \frac{1}{\sin^{2}(\frac{k\pi}{N})} = \frac{N^{2}-1}{3}, \qquad 
\sum_{k=1}^{N-1} \frac{1}{\sin^{4}(\frac{k\pi}{N})} = \frac{(N^{2}-1)(N^{2}+11)}{45 N}. $$
\end{corollary}

{\bf Proof.}
Since ${\bf A}$ in (\ref{A1}) resp. ${\bf A}^{2}$ are  circulant matrices, they can be diagonalized by the Fourier matrix ${\bf F}_{N} = \frac{1}{\sqrt{N}} \left( \omega_{N}^{jk} \right)_{j,k=0}^{N-1}$, where 
$\omega_{N} := e^{-2\pi i/N}$.
We obtain $\widehat{\bf A} =  {\bf F}_{N} {\bf A} {\bf F}_{N}^{*} = {\rm diag} \, 
\left(a(\omega_{N}^{k}) \right)_{k=0}^{N-1}$ with the characteristic polynomial $a(z):= -2 + z + z^{N-1}$ of ${\bf A}$, i.e.,
$$ a(\omega_{N}^{k}) = -2 + e^{-2 \pi i k /N} + e^{2 \pi i k /N} = -2+ 2 \cos \left( \frac{2\pi k}{N} \right) = - 4 \sin^{2}\left(\frac{\pi k }{N} \right) , \quad k=0, \ldots , N-1. $$
Observe that $a(\omega_{N}^{k})$ are the eigenvalues of ${\bf A}$ resp.~$\widehat{\bf A}$.
Considering the generalized Moore-Penrose inverse, we obtain the circulant matrix ${\bf A}^{+} = {\bf F}_{N}^{*} \widehat{\bf A}^{+} {\bf F}_{N}$ with 
$$ \widehat{\bf A}^{+} := {\rm diag} \, \left( 0, \frac{-1}{4 \sin^{2}(\frac{\pi}{N})}, \frac{-1}{4 \sin^{2}(\frac{2\pi}{N})}, \ldots , \frac{-1}{ 4 \sin^{2}(\frac{{(N-1)}\pi}{N})} \right).$$ 
Analogously, the Moore-Penrose inverse of ${\bf A}^{2}$ can be written as $({\bf A}^{2})^{+} = {\bf F}_{N}^{*} (\widehat{\bf A}^{+})^{2} {\bf F}_{N}$. 
A comparison of the entries $a_{j}^{+}$ (resp.\ $c_{j}^{+}$) in the first column of of ${\bf F}_{N}^{*} \widehat{\bf A}^{+} {\bf F}_{N}$, (resp.\ ${\bf F}_{N}^{*} (\widehat{\bf A}^{+})^{2} {\bf F}_{N}$),
$$ a_j^{+} = \frac{1}{N} \sum_{k=1}^{N-1} \frac{-1}{4 \sin^{2}(\frac{\pi k}{N})} \omega_N^{-jk} \omega_N^{k0}
= \frac{1}{N}\sum_{k=1}^{N-1} -\frac{\cos{\frac{2\pi kj}{N}}+i\sin{\frac{2\pi kj}{N}}}{4\sin^{2}(\frac{k\pi}{N})}
$$
resp. 
$$ c_j^{+} =  \frac{1}{N}\sum_{k=1}^{N-1} -\frac{\cos{\frac{2\pi kj}{N}}+i\sin{\frac{2\pi kj}{N}}}{16\sin^{4}(\frac{k\pi}{N})}, $$
with the formulas found for ${\bf A}^{+}$ (resp.\  (${\bf A}^{2})^{+}$) in Theorems \ref{theo2} and \ref{theobi} yields the assertions.
\hfill $\Box$  


\medskip

Let us now consider the case of reflecting boundary conditions and derive the Moore-Penrose inverse for the matrix ${\bf B}$ in (\ref{A2}) and for $-{\bf B}^{2}$.
\begin{theorem}\label{theo3}
The Moore-Penrose inverse of the matrix ${\bf B}$ in {\rm (\ref{A2})} occurring by discretization of the Laplacian in the case of homogeneous Neumann boundary conditions is a symmetric matrix of the form ${\bf B}^{+} = (b_{j,k}^{+})_{j,k=0}^{N-1}$ with 
$$ b_{j,k}^{+} = - \frac{(N-1)(2N-1)}{6N}  + j - \frac{j(j+1)}{2N} - \frac{k(k+1)}{2N}$$
for $j \ge k$.
\end{theorem}

{\bf Proof}.
1. Let now ${\bf A} = {\bf A}_{2N} = {\rm circ}\,  (-2,1, \ldots , 1)^{T} \in {\IR}^{2N \times 2N}$ be the circulant matrix as in (\ref{A1}) but of order $2N$.
  We  consider the four $N \times N$ blocks of ${\bf A}_{2N}$,
$$ {\bf A}_{2N} = \left[ \begin{array}{cc} {\bf A}_{0} & {\bf A}_{1}^{T} \\ {\bf A}_{1} & {\bf A}_{0} \end{array} \right], $$
where ${\bf A}_{0}$ and ${\bf A}_{1}$  are symmetric Toeplitz matrices of size $N \times N$.
In particular, ${\bf A}_{1} = {\bf A}_{1}^{T}$ has only two nonzero entries. Further, we denote by  
${\bf A}^{+}_{2N} = {\rm circ} \, (a_{0}^{+}, a_{1}^{+}, \ldots , a_{2N-1}^{+})$
the Moore-Penrose inverse of ${\bf A}_{2N}$, where $a_{k}^{+} = \frac{1}{24 N} (1 - 4 N^{2}) + \frac{1}{4N} k(2N-k)$ by Theorem \ref{theo2}. Let the $N \times N$ blocks of ${\bf A}_{2N}^{+}$ be denoted by 
$${\bf A}^{+}_{2N}=
\left[ \begin{array}{cc} {\bf A}_{0}^{+} & {\bf A}_{1}^{+T} \\ {\bf A}_{1}^{+} & {\bf A}_{0}^{+}\end{array} \right]. $$
 Observe that also 
 ${\bf A}_{0}^{+}$ and ${\bf A}_{1}^{+}$ are Toeplitz matrices. Further, let $\tilde{\bf I}$ denote the counter identity of size $N \times N$, i.e.,
$$ \tilde{\bf I} = \left( \begin{array}{cccc}   &  &  & 1 \\
&  &  1 &  \\ & \varddots & &  \\ 1 &  & & \end{array} \right). $$
Then we simply observe that the matrix ${\bf B}$ in (\ref{A2}) can be written as 
$$ {\bf B} = {\bf A}_{0} + \tilde{\bf I} {\bf A}_{1}. $$
We show now that  ${\bf B}^{+} := {\bf A}_{0}^{+} + \tilde{\bf I} {\bf A}_{1}^{+}$ is the Moore-Penrose inverse of ${\bf B}$, where  ${\bf A}_{0}^{+}$  and ${\bf A}_{1}^{+}$ are the partial matrices of ${\bf A}_{2N}^{+}$ as given above. 
We know that
$$
{\bf A}_{2N} {\bf A}^{+}_{2N} {\bf A}_{2N} = \left[ \begin{array}{cc} {\bf A}_{0} & {\bf A}_{1}^{T} \\ {\bf A}_{1} & {\bf A}_{0} \end{array} \right]    \left[ \begin{array}{cc} {\bf A}_{0}^{+} & {\bf A}_{1}^{+T} \\ {\bf A}_{1}^{+} & {\bf A}_{0}^{+}\end{array} \right]        \left[ \begin{array}{cc} {\bf A}_{0} & {\bf A}_{1}^{T} \\ {\bf A}_{1} & {\bf A}_{0} \end{array} \right] = {\bf A}_{2N} = \left[ \begin{array}{cc} {\bf A}_{0} & {\bf A}_{1}^{T} \\ {\bf A}_{1} & {\bf A}_{0} \end{array} \right]$$
which leads to
\begin{eqnarray*}
{\bf A}_{0} {\bf A}_{0}^{+} {\bf A}_{0} + {\bf A}_{1}^{T} {\bf A}_{1}^{+} {\bf A}_{0} + {\bf A}_{0} ({\bf A}_{1}^{+})^{T} {\bf A}_{1} + {\bf A}_{1}^{T} {\bf A}_{0}^{+} {\bf A}_{1} &=& {\bf A}_{0}, \\
{\bf A}_{1} {\bf A}_{0}^{+} {\bf A}_{0} + {\bf A}_{0} {\bf A}_{1}^{+} {\bf A}_{0} + {\bf A}_{1} ({\bf A}_{1}^{+})^{T} {\bf A}_{1} + {\bf A}_{0} {\bf A}_{0}^{+} {\bf A}_{1} &=& {\bf A}_{1}.
\end{eqnarray*}
Using these equalities, and taking into account the symmetry of ${\bf A}_{0}$, ${\bf A}_{1}$ and ${\bf A}_{0}^{+}$ as well as that $\tilde{\bf I} {\bf T} \tilde{\bf I} = {\bf T}^{T}$ for any Toeplitz matrix ${\bf T}$, a short computation gives
\begin{eqnarray*}
 {\bf B} {\bf B}^{+} {\bf B} &=& ({\bf A}_{0} + \tilde{\bf I} {\bf A}_{1}) ({\bf A}_{0}^{+} + \tilde{\bf I} {\bf A}_{1}^{+}) ({\bf A}_{0} + \tilde{\bf I} {\bf A}_{1}) \\
 &=& {\bf A}_{0} {\bf A}_{0}^{+} {\bf A}_{0} + {\bf A}_{0} {\bf A}_{0}^{+} \tilde{\bf I} {\bf A}_{1} + {\bf A}_{0} \tilde{\bf I} {\bf A}_{1}^{+} {\bf A}_{0} +
 {\bf A}_{0} \tilde{\bf I} {\bf A}_{1}^{+} \tilde{\bf I} {\bf A}_{1} + \tilde{\bf I} {\bf A}_{1} {\bf A}_{0}^{+} {\bf A}_{0} + \tilde{\bf I} {\bf A}_{1} {\bf A}_{0}^{+} \tilde{\bf I} {\bf A}_{1}\\
 & &  + \tilde{\bf I} {\bf A}_{1} \tilde{\bf I} {\bf A}_{1}^{+} {\bf A}_{0} + \tilde{\bf I} {\bf A}_{1} \tilde{\bf I} {\bf A}_{1}^{+} \tilde{\bf I} {\bf A}_{1} \\
 &=& {\bf A}_{0} {\bf A}_{0}^{+} {\bf A}_{0} + \tilde{\bf I}  {\bf A}_{0}^{T} ({\bf A}_{0}^{+})^{T} {\bf A}_{1}^{T} + \tilde{\bf I} {\bf A}_{0}^{T} {\bf A}_{1}^{+} {\bf A}_{0}
 + {\bf A}_{0} ({\bf A}_{1}^{+})^{T} {\bf A}_{1} + \tilde{\bf I} {\bf A}_{1} {\bf A}_{0}^{+} {\bf A}_{0}  \\
 & & + {\bf A}_{1}^{T} ({\bf A}_{0}^{+})^{T} {\bf A}_{1} + {\bf A}_{1}^{T} {\bf A}_{1}^{+} {\bf A}_{0} + \tilde{\bf I} {\bf A}_{1} ({\bf A}_{1}^{+})^{T} {\bf A}_{1} \\
 &=&
({\bf A}_{0} + \tilde{\bf I} {\bf A}_{1}) = {\bf B}
\end{eqnarray*}
and analogously ${\bf B}^{+} {\bf B} {\bf B}^{+} = {\bf B}^{+}$. 
Further, using the symmetry of ${\bf A}_{2N} {\bf A}_{2N}^{+}$ we observe that $({\bf A}_{0} {\bf A}_{0}^{+} + {\bf A}_{1}^{T} {\bf A}_{1}^{+})$ is symmetric 
yielding
\begin{eqnarray*}
{\bf B} {\bf B}^{+} &=& ({\bf A}_{0} + \tilde{\bf I} {\bf A}_{1}) ({\bf A}_{0}^{+} + \tilde{\bf I} {\bf A}_{1}^{+}) \\
&=& {\bf A}_{0} {\bf A}_{0}^{+} + \tilde{\bf I} {\bf A}_{1} \tilde{\bf I} {\bf A}_{1}^{+} + {\bf A}_{0} \tilde{\bf I} {\bf A}_{1}^{+} + \tilde{\bf I} {\bf A}_{1} {\bf A}_{0}^{+} \\
&=& ({\bf A}_{0} {\bf A}_{0}^{+} +  {\bf A}_{1}^{T} {\bf A}_{1}^{+}) + \tilde{\bf I} ({\bf A}_{0}^{T} {\bf A}_{1}^{+} + {\bf A}_{1} {\bf A}_{0}^{+}) \\
&=& ({\bf A}_{0} {\bf A}_{0}^{+} +  {\bf A}_{1}^{T} {\bf A}_{1}^{+})^{T} +   ({\bf A}_{0}^{T} {\bf A}_{1}^{+} + {\bf A}_{1} {\bf A}_{0}^{+})^{T} \tilde{\bf I} \\
&=& ({\bf A}_{0}^{+})^{T} {\bf A}_{0}^{T} +  ({\bf A}_{1}^{+})^{T} {\bf A}_{1} +  ({\bf A}_{1}^{+})^{T} {\bf A}_{0} \tilde{\bf I} +  ({\bf A}_{0}^{+})^{T} {\bf A}_{1}^{T} \tilde{\bf I} \\
&=& ({\bf A}_{0}^{+} + \tilde{\bf I} {\bf A}_{1}^{+})^{T} ({\bf A}_{0} + \tilde{\bf I} {\bf A}_{1})^{T} = ({\bf B} {\bf B}^{+})^{T}
\end{eqnarray*}
and analogously ${\bf B}^{+} {\bf B} =
({\bf B}^{+} {\bf B})^{T}$.

2. Using the formulas for the $(2N \times 2N)$ matrix ${\bf A}_{2N}^{+}$ from Theorem \ref{theo2}  we now obtain ${\bf B}^{+} = (b_{j,k}^{+})_{j,k=0}^{N-1}$  from  ${\bf A}_{2N}^{+} = (a_{j,k}^{+})_{j,k=0}^{2N-1} = (a_{(j-k) {\rm mod} \, 2N}^{+})_{j,k=0}^{2N-1}$ with
$$ b_{j,k}^{+} = a_{j,k}^{+} + a_{2N-1-j,k}^{+} = a^{+}_{(j-k) {\rm mod}\,  2N} +
a^{+}_{(2N-j-k-1) {\rm mod} \,  2N}. $$
For $j \ge k$ this gives 
\begin{eqnarray*}
b_{j,k}^{+} &=& \frac{2(1-(2N)^{2})}{24 N} + \frac{(j-k)(2N-j+k)}{4N} + \frac{(2N-j-k-1)(j+k+1)}{4N} \\
&=& - \frac{(N-1)(2N-1)}{6N} + j - \frac{j(j+1)}{2N} - \frac{k(k+1)}{2N}.
\end{eqnarray*}
\null \hfill $\Box$

\begin{theorem}\label{theobi1}
The Moore-Penrose inverse for the discrete biharmonic operator ${\bf L}=-{\bf B}^{2}$ with ${\bf B}$ in {\rm (\ref{A2})} 
is a symmetric matrix of the form ${\bf L}^{+}= -({\bf B}^{2})^{+} = (d_{j,k}^{+})_{j,k=0}^{N-1}$ with 
\begin{eqnarray*}
 d_{j,k}^{+} &=& 
  - \frac{(N^{2}-1)(4N^{2}-1)}{180N}  
  + \frac{\left[j(j+1)+k(k+1)\right]^{2}}{24N} + \frac{j(j+1)k(k+1)}{6N} \\
 & &  + \frac{N}{6} \left[(j(j+1)+k(k+1) \right] - \frac{(2j+1)}{12} \left[j(j+1)+3k(k+1)\right]
 \end{eqnarray*}
for $j \ge k$.
\end{theorem}

{\bf Proof}. 
From Theorem \ref{theo3} it follows that $b_{j,k}^{+} = \beta_{j,k} + j$ for $j \ge k$ and $b_{j,k}^{+} = \beta_{j,k} + k$ for $j < k$, where
$$ \beta_{j,k} := - \frac{(N-1)(2N-1)}{6N} - \frac{j(j+1)}{2N} -\frac{k(k+1)}{2N}. $$
Observing that $({\bf B}^{2})^{+} = ({\bf B}^{+})^{2}$, we obtain for $j \ge k$
$$
 -d_{j,k}^{+} = \sum_{\ell=0}^{k} (\beta_{j,\ell} + j)(\beta_{\ell,k} + k) + \sum_{\ell=k+1}^{j} (\beta_{j,\ell} + j)(\beta_{\ell,k} + \ell)
+ \sum_{\ell=j+1}^{N-1} (\beta_{j,\ell} + \ell)(\beta_{\ell,k} + \ell). $$
Inserting the value $\beta_{j,\ell}$, some algebraic evaluations give the desired result.
\hfill $\Box$

\section{Sparse signal approximation by regularized diffusion inpainting}
\label{sec4}

Using the above explicit representations of the Moore-Penrose inverses of 
difference matrices, we want to derive a new algorithm for sparse signal 
approximation based on a regularization of a discrete diffusion inpainting 
method.

\subsection{Regularized diffusion inpainting}

We consider now the following approach for sparse signal representation.
For a finite signal ${\bf f} = (f_{j})_{j=0}^{N-1} \in {\IR}^{N}$, we want to find a good approximation ${\bf u} = (u_{j})_{j=0}^{N-1} \in {\IR}^{N}$ using only a small amount of values $\widetilde{f}_{k}$, with $k \in \Gamma \subset \Omega:= \{ 0, \ldots , N-1 \}$.
Here, $\Gamma$ denotes  an index set with a small number $|\Gamma|$ of entries,
i.e. $|\Gamma| \ll N$. 
In order to construct ${\bf u}$, we apply an inpainting procedure.  We want to approximate ${\bf f}$ by finding a solution of 
\begin{eqnarray} \label{i1}
( {\bf L} {\bf u})_{k} &=& 0 \quad {\rm for} \; k \in \{ 0, \ldots , N-1 \} \setminus \Gamma,\\
\label{i2}  u_{k} &=& \widetilde{f}_{k} \quad {\rm for} \; k \in \Gamma,
\end{eqnarray}
where ${\bf L}$ denotes a symmetric difference matrix satisfying the conditions (i) and (ii). Here, ${\bf L}$  can be seen as a smoothing operator, as e.g. the discretization of the second derivative (Laplacian) or the fourth derivative (biharmonic operator) with periodic or reflecting boundary conditions. 

Let the vector ${\bf c} := (c_{0}, \ldots , c_{N-1})^{T}$ be given by
$$c_{i} := \left\{\begin{array}{ll}  1 & i \in \Gamma, \\
0 & i \in \Omega \setminus \Gamma. \\  \end{array} \right. $$
Further, let ${\bf C} := {\rm diag} \, ({\bf c})$ be the diagonal matrix determined by ${\bf c}$, and let ${\bf I}_{N}$ denote the identity matrix of order $N$.

No we can reformulate the problem as follows. Find an index set $\Gamma$ (or equivalently an index vector ${\bf c})$ with a fixed number of values $|\Gamma | \ll N$, as well as a sparse vector $\widetilde{\bf f} = (\widetilde{f}_{0}, \ldots , \widetilde{f}_{N-1})^{T}$  with $\widetilde{f}_{k} =0$ for $k \in \Omega \setminus \Gamma$,
such that ${\bf u}$ can be recovered from $\widetilde{\bf f}$
as
\begin{equation} \label{e1}
{\bf C} ( {\bf u} - \widetilde{\bf f})  - ({\bf I}_{N} - {\bf C}) {\bf L} {\bf u} = {\bf 0}
\end{equation}
and $\| {\bf u} - {\bf f} \|_{2}$ is (in some sense) minimal.

In the two-dimensional case with homogeneous Neumann boundary conditions, Mainberger et al.~\cite{MHW12} considered a specific greedy algorithm: To
determine a suitable set $\Gamma$, they used a probabilistic sparsification followed by a nonlocal pixel exchange. Afterwards, the tonal (i.e. greyvalue) data 
$\widetilde{f}(k)$ for $k \in \Gamma$ are optimized with a least squares approach.
Here we propose a different strategy that simultaneously optimizes spatial and tonal data. Equation (\ref{e1}) yields
$$ {\bf L} {\bf u} = {\bf C} ({\bf u} - \tilde{\bf f} + {\bf L} {\bf u}). $$
Since  ${\bf C}$ is a sparse diagonal matrix with only $|\Gamma|$ nonzero diagonal entries, it follows that ${\bf g}:= {\bf C} ({\bf u} - \tilde{\bf f} + {\bf L} {\bf u})$ is a sparse vector with $|\Gamma|$ nonzero entries. Therefore, we 
reformulate our problem: We aim to solve 
$$ {\bf L} {\bf u} = {\bf g}, $$
where the sparse vector ${\bf g}$ needs to be chosen in a way  such that the error $\| {\bf u} - {\bf f} \|_{2}$ is small.
Since ${\bf L}$ is not invertible we apply the Moore-Penrose inverse of ${\bf L}$ to obtain 
\begin{equation} \label{L+} {\bf u} = {\bf L}^{+} {\bf g}.
\end{equation}
From  ${\bf 1}^{T} {\bf L} = {\bf 0}^{T}$, we find that  ${\bf 1}^{T} {\bf L} {\bf u} = {\bf 1}^{T} {\bf g} = 0$. This solvability condition is satisfied  if ${\bf 1}^{T} {\bf C}({\bf u} - \tilde{\bf f} - {\bf L} {\bf u})=0$ is true.
Further, we also have ${\bf 1}^{T} {\bf L}^{+} ={\bf 0}^{T}$, i.e., each solution ${\bf u} = {\bf L}^{+} {\bf g}$ in (\ref{L+}) has the mean value $0$. Since we are interested in a sparse approximation ${\bf u}$ of ${\bf f}$, we therefore also assume that also ${\bf f}$ has mean value $0$. In practice, we just replace ${\bf f}$ by  ${\bf f}^{0} = {\bf f} - \bar{f} {\bf 1}$ with 
 $$ \bar{f}:=\frac{1}{N}{\bf 1}^{T} {\bf f} = \frac{1}{N}\sum_{j=0}^{N-1} {f}_{j}  $$
 and store $\bar{f}$ elsewhere.
 
 
\medskip

\noindent 
{\bf Remark.} 
The regularization approach described above to solve the inpainting problem can also be interpreted  as a solution using the concept of 
discrete Green's functions; see \cite{HPW15} for the two-dimensional case.


\subsection{Orthogonal matching pursuit approach}

With the regularization approach (\ref{L+}), the original task of finding  a suitable index set $\Gamma$ and the corresponding optimal  values $\tilde{f}_{k}$, $k \in \Gamma$ can be rewritten  as the following optimization problem:
Find  ${\bf g} \in {\IR}^{N}$ such that 
$$ \| {\bf L}^{+} {\bf g} - {\bf f} \|_{2} \to \min \quad \hbox{subject to} \quad \| {\bf g} \|_{0} = |\Gamma|, $$
where $\| {\bf g} \|_{0}$ denotes the number of nonzero coefficients in the vector ${\bf g}$.
Let us denote  the columns of the matrix $ {\bf L}^{+}$  by $\textbf{a}_0,\ldots,\textbf{a}_{N-1}$. 
We now rewrite the optimization problem in the form,
  $$  \|{\bf \nufett } \|_{2} \to \min \quad   \hbox{subject to} \quad  {\bf L}^{+} {\bf g} + \nufett = {\bf f}, \; \| {\bf g} \|_{0} =  |\Gamma|.  $$
Equivalently, with $\Gamma$ denoting the index set of the nonzero entries of the sparse vector ${\bf g}$ and $|\Gamma| = \| {\bf g} \|_{0} \ll N$, we have ${\bf f}$ in the form
\begin{eqnarray*}
{\textbf{f}}=\sum\limits_{k \in {\Gamma}} g_{k} \, \textbf{a}_k + \nufett,
\end{eqnarray*}
where ${\textbf{f}}$ is approximated by a sparse linear combination of the columns $\textbf{a}_k$, $k \in \Gamma$. 
We apply now the orthogonal matching pursuit as a greedy algorithm for finding a suitable subset ${\Gamma}$ as well as the sparse vector ${\bf g}$ that determines ${\bf u}$.
Here,  $\{ {\bf a}_{0}, \ldots , {\bf a}_{N-1} \}$ is the dictionary for our OMP method in the Hilbert space ${\IR}^{N}$.
This algorithm is known to work very efficiently also when the original signal ${\bf f}$ is not exactly sparse in the considered dictionary, see e.g. \cite{TG07}.
\medskip

The OMP works as follows.
In a first step, we determine the index $k_1\in\{0,\ldots, N-1\}$ 
such that the column $\textbf{a}_{k_1}$ correlates most strongly with ${\textbf{f}}$, i.e.
\begin{eqnarray*}
k_1=\rm{arg}\max\limits_{k=0,\ldots,N-1} \frac{|\langle {\textbf{f}},\textbf{a}_k\rangle| }{\langle {\textbf{a}}_{k},\textbf{a}_k\rangle},
\end{eqnarray*}
where $\langle {\textbf{f}},\textbf{a}_k\rangle= {\bf f}^{T} {\bf a}_{k}$ is the standard scalar product of the two vectors ${\textbf{f}}$ and $\textbf{a}_k$.\\
In the next step, we determine the coefficient $g_{k_{1}}$ such that the Euclidean norm 
$\|{\textbf{f}}- g_{k_{1}}\cdot\textbf{a}_{k_1}\|_2$ is minimal, i.e. 
$g_{k_{1}}=\langle {\textbf{f}},\textbf{a}_{k_1}\rangle / \| \textbf{a}_{k_{1}}\|_{2}^{2}$, 
where $\| \textbf{a}_{k_{1}}\|_{2}$ denotes the Euclidean norm of $\textbf{a}_{k_{1}}$.

Now we consider the residuum $\textbf{r}_1= {\textbf{f}}-g_{k_{1}}\textbf{a}_{k_1}$ and proceed with the next iteration step, where ${\bf f}$ is replaced by ${\bf r}_{1}$. For all further iteration steps we add an update  of the coefficients. The algorithm can be summarized as follows.

\medskip

\begin{algorithm} (OMP) \label{alg1} \\
{\bf Input:} Dictionary $\{ {\bf a}_{0}, \ldots, {\bf a}_{N-1} \}$ of ${\IR}^{N}$, ${\bf r}_0={\bf f} \in {\IR}^{N}$ (with ${\bf 1}^{T} {\bf f} = 0$), ${\bf g}={\bf 0}$, $|\Gamma| \ll N$.
\medskip

\noindent
For $j=1, \ldots , |\Gamma|$ do 
\begin{enumerate}
\item Determine an optimal index $k_j$ such that ${\bf a}_{k_j}$ correlates most strongly with the residuum ${\bf r}_{j-1}$, i.e.
\begin{eqnarray*}
k_j={\rm arg }\max\limits_{k=0,\ldots,N-1}\frac{|\langle {\bf r}_{j-1},{\bf a}_k\rangle|}{\langle {\bf a}_{k},{\bf a}_k\rangle}.
\end{eqnarray*}
\item Update the coefficients $g_{k_{1}},\ldots,g_{k_{j}}$ such that $\|{\bf f}-\sum_{i=1}^j g_{k_{i}}\, {\bf a}_{k_i}\|_2$ is minimal
 and set ${\bf r}_j={\bf f}-\sum_{i=1}^j g_{k_{i}}{\bf a}_{k_i}$.
\end{enumerate}
end (do)
\smallskip

\noindent
{\bf Output:} ${\bf g}^{\Gamma} = (g_{k_{i}})_{i=1}^{|\Gamma|}$
\end{algorithm}

\noindent

As proposed in the algorithm, for a given number $|\Gamma| \ll N$, we may just take $|\Gamma|$ iteration steps. For improving the method, one may also take some more steps, and afterwards decide by a shrinkage procedure, which indices are kept to determine $\Gamma$.

\medskip

\noindent
{\bf Remarks.} 
1. For numerical purposes, we normalize the columns ${\bf a}_{0}, \ldots , {\bf a}_{N-1}$ in a preprocessing step before starting the OMP algorithm. Using ${\bf L}^{+}$ based on the Laplace operator or the biharmonic operator with periodic boundary conditions, see Theorems \ref{theo2} and \ref{theobi}, all columns have the same Euclidean norm, and we obtain for example
$$ \sum_{k=0}^{N-1} (a_{k}^{+})^{2} = \sum_{k=0}^{N-1} \left( \frac{1-N^{2}}{12N} + \frac{k(N-k)}{2N} \right)^{2} = \frac{(N^{2}+11)(N^{2} -1)}{720 N} $$
for the Moore-Penrose inverse in Theorem \ref{theo2}.
Note that the columns occurring in the Moore-Penrose inverses in Theorems \ref{theo3} and \ref{theobi1} for reflecting boundary conditions do not have the same Euclidean norm, such that the normalization in the OMP algorithm is crucial.

2. Observe that the considered regularization does not longer exactly solve the original inpainting problem (\ref{i1})--(\ref{i2}) resp.\ (\ref{e1}),
since we cannot ensure the solvability condition ${\bf 1}^{T} {\bf g}  = {\bf 0} $.
One way to enforce this condition is by incorporating it as additional constraint, see Section 4.

\subsection{Relation to spline interpolation with arbitrary knots}

Finally, we want to point out the close relationship between the proposed inpainting algorithm considered above and the nonlinear spline approximation problem with arbitrary knots in order to find a sparse approximation of ${\bf f}$.
For that purpose, we consider the inpainting problem (\ref{e1}) for the special case ${\bf L} = {\bf B}$ in (\ref{A2}) corresponding to the Laplace operator with reflecting (homogeneous Neumann) boundary conditions,
\begin{equation}\label{ef}
 {\bf C} {\bf u} - ({\bf I}_{N} - {\bf C}){\bf B} {\bf u} = {\bf C} \widetilde{\bf f}.
 \end{equation}
Assume that the index set of $1$-entries of the index vector ${\bf c}$ is $\Gamma = \{ \ell_{1}, \ldots , \ell_{L} \}$, where $0 \le \ell_{1} < \ell_{2} < \ldots  < \ell_{L} \le N-1$.
In this case, the matrix ${\bf M} := {\bf C} - ({\bf I}_{N} - {\bf C}) {\bf B}$ is an invertible block diagonal matrix with $L+1$ blocks ${\bf M}_{j}$, where the inner blocks are of the form
$$ {\bf M}_{j} = \left( \begin{array}{ccccc} 1 & 0 & \ldots  &  & 0 \\
1 & -2 & 1 &  & 0 \\
 & \ddots & \ddots & \ddots  & \\
 0 &  & 1 & -2 & 1 \\
 0 & \ldots & 0 & 0 &  1 \end{array} \right) \in {\IR}^{(\ell_{j+1} - \ell_{j}+1) \times (\ell_{j+1} - \ell_{j}+1)}, \qquad \ell=1, \ldots , L-1, $$
 and the boundary blocks
 $$ {\bf M}_{0} = \left( \begin{array}{ccccc} -1 & 1 & 0 &\ldots    & 0 \\
1 & -2 & 1 &  & 0 \\
 & \ddots & \ddots & \ddots  & \\
 0 &  & 1 & -2 & 1 \\
 0 & \ldots & 0 & 0 &  1 \end{array} \right) 
\qquad  {\bf M}_{L} = \left( \begin{array}{ccccc} 1 & 0 & \ldots  &  & 0 \\
1 & -2 & 1 &  & 0 \\
 & \ddots & \ddots & \ddots  & \\
 0 &  & 1 & -2 & 1 \\
 0 & \ldots & 0 & 1 &  -1 \end{array} \right) 
 $$
are of order $\ell_{1}+1$ resp. $N - \ell_{L}$.
Observe that the boundary blocks can degenerate to ${\bf M}_{0} = 1$ for $\ell_{1} =0$ and ${\bf M}_{0} = \left( \begin{array}{cc} -1 & 1 \\ 0 & 1 \end{array} \right)$ for $\ell_{1} = 1$, and analogously, ${\bf M}_{L} = 1$ for $\ell_{L} =N-1$ and ${\bf M}_{L} = \left( \begin{array}{cc} 1 & 0 \\ 1 & -1 \end{array} \right)$ for $\ell_{L} = N-2$.
The inner blocks degenerate to
 ${\bf M}_{j} = {\bf I}_{2}$ for $\ell_{j+1} = \ell_{j}+1$.
Now the system ${\bf M} {\bf u} = {\bf C} \widetilde{\bf f}$ in (\ref{ef}) can be solved by solving the partial systems 
$${\bf M}_{j} {\bf u}^{(j)} = {\bf C}_{j} \widetilde{\bf f}^{(j)}, $$ 
where 
 \begin{eqnarray*}
 {\bf C}_{j} &:=& {\rm diag} (1, 0, \ldots , 0, 1) \in {\IR}^{(\ell_{j+1} - \ell_{j}+1) \times (\ell_{j+1} - \ell_{j}+1)} \qquad  \ell=1, \ldots , L-1, \\
 {\bf C}_{1} &:=& {\rm diag} (0, \ldots , 0, 1) \in {\IR}^{(\ell_{1}+1) \times (\ell_{1}+1)}, \\
  {\bf C}_{L} &:=& {\rm diag} (1, 0, \ldots , 0) \in {\IR}^{(N - \ell_{L}) \times (N - \ell_{L})}, \\
{\bf u}^{(j)} &:=& (u_{\ell_{j}}, \ldots u_{\ell_{j+1}})^{T}, \qquad j=0, \ldots , L\\
 \widetilde{\bf f}^{(j)} &:=& (\widetilde{f}_{\ell_{j}}, 0 , \ldots , 0,  \widetilde{f}_{\ell_{j+1}})^{T} \qquad j=0, \ldots , L.
\end{eqnarray*}
 For the inner blocks we observe that ${\bf C}_{j} \widetilde{\bf f}^{(j)} = (\widetilde{f}_{\ell_{j}}, 0  , \ldots , 0, \widetilde{f}_{\ell_{j+1}})^{T}$, and it can be simply verified 
 that 
 $$ u_{\ell_{j}+k} = \frac{(n_{j}-k) \widetilde{f}_{\ell_{j}} + k \widetilde{f}_{\ell_{j+1}}}{n_{j}}, \qquad n_{j}:= \ell_{j+1} - \ell_{j}, \; k=0, \ldots , n_{j}. $$
For the boundary blocks, the constant vectors 
$$ {\bf u}^{(0)} = \widetilde{f}_{\ell_{1}} \, {\bf 1}_{\ell_{1}}, \qquad {\bf u}^{(L)} = \widetilde{f}_{\ell_{L}} \, {\bf 1}_{N-\ell_{1}} $$
solve the partial systems.
As a global solution, we hence  obtain a discrete linear spline approximation ${\bf u}$ of ${\bf f}$ with knots given by the index set $\Gamma$ resp.\ ${\bf c}$.
Thus, the inpainting problem considered in (\ref{ef}) is equivalent to the nonlinear approximation problem to find optimal knots $\ell_{j}$ (out of the set $\{ 0, \ldots , N-1 \})$ and corresponding 
values $\widetilde{f}_{\ell_{j}}$, such that the obtained discrete piecewise linear spline vector 
$$ {\bf u} = ( ({\bf u}^{(0)})^{T}, \ldots , ({\bf u}^{(L)})^{T} )^{T} \in {\IR}^{N} $$
minimizes the norm $\| {\bf f} - {\bf u} \|_{2}$.

As  we will see in the next section, our regularization approach using the OMP method can be thus also seen as a method to find an solution for the discrete spline approximation problem. 
\medskip

\noindent
{\bf Remarks.}\\
1. Analogously, it can be shown that the inpainting algorithm with the operator $-{\bf B}^{2}$ corresponding to the biharmonic operator with reflecting boundary conditions is related to  nonlinear discrete spline approximation with cubic splines with arbitrary knots.

2. Sparse data approximation using splines with free knots has been already extensively studied in the literature, see e.g. {\cite{J78,LW91,SS95,B04}} and references therein. Splines with free knots have been shown to perform significantly better than splines with fixed knots. For a discussion in the context of signal
approximation we refer to \cite{Ho15}.
The question to find the optimal knot sequence leads to a global optimization problem that is nonlinear and non-convex.
Therefore, besides expensive global optimization strategies, there exist also heuristic strategies, where one starts either with  only a few knots  and iteratively adds further knots to improve the quality of the approximating spline, or conversely, starts with a large number of knots and iteratively removes knots that are less important for the approximation error.
Our approach belongs to the first group.
 

\section{Numerical examples}
We apply now the regularization of the inpainting approach  proposed in Section 3 to sparse approximation of vectors. In particular, we want to use the Moore-Penrose inverses of
the discrete Laplacian and the discrete biharmonic operator  with either  periodic or reflecting boundary conditions that have been explicitly computed in Section 2. 
Our previous considerations imply the following simple algorithm.

\begin{algorithm} (Sparse approximation)\\
{\bf Input:} $\{  {\bf a}_{0}, \ldots ,  {\bf a}_{N-1} \}$ columns of the matrix  ${\bf L}^{+}$,\\
\null \hspace{13mm} ${\bf f} \in {\IR}^{N}$, $|\Gamma| +1$ number of values used for sparse approximation 
\smallskip

\begin{enumerate}
\item Compute $\bar{f}:=\frac{1}{N}{\bf 1}^{T} {\bf f}$ and ${\bf f}^{0} := {\bf f} - \bar{f}{\bf 1}$.
\item Apply the OMP-algorithm \ref{alg1} to ${\bf f}^{0}$ with $|\Gamma|$ steps to obtain an approximation $\sum\limits_{i=1}^{|\Gamma|} g_{k_{i}} \, {\bf a}_{k_{i}}$ of ${\bf f}^{0}$.
\end{enumerate}

\noindent
{\bf Output:} $\bar{f}$, ${\bf g}^{\Gamma} =(g_{k_{i}})_{i=1}^{|\Gamma|}$ determining a sparse approximation of ${\bf f}$ of the form 
$$ {\bf u} = \bar{f} {\bf 1} + \sum_{i=1}^{|\Gamma|} g_{k_{i}} \, {\bf a}_{k_{i}}. $$
\end{algorithm}

We have applied this algorithm  using the columns of the Moore-Penrose inverse of ${\bf L}$ for 
${\bf L} \in \{ {\bf A},  \, {\bf B}, \, - {\bf A}^{2}, \, - {\bf B}^{2} \}$ as computed on Section 2. Starting with a vector ${\bf f}$ of length $N=256$, we spend $|\Gamma| = 13$ values for the approximation of ${\bf f}$, i.e., about 5 \% of the original vector size.
In Figure 1, we present  the approximation results together with their error. 
As is common in signal processing, we measure the error in terms of its
peak signal-to-noise ratio (PSNR). For a maximum possible signal value of $M$ 
(in our experiments: $M=1$), the PSNR between two signals 
${\bf u}=(u_i)_{i=0}^{N-1}$ and ${\bf f}=(f_i)_{i=0}^{N-1}$ is given by
\begin{equation}
 {\rm PSNR}\,({\bf u},{\bf f}) \;:=\;
   10 \, \log_{10} \, \left( \frac{M^2}{{\rm MSE}\,({\bf u},{\bf f})} \right)
\end{equation}
where MSE denotes the mean squared error
\begin{equation}
 {\rm MSE}\,({\bf u},{\bf f}) = \frac{1}{N} \sum_{i=0}^{N-1} (u_i - f_i)^2.
\end{equation}
Note that higher PSNR values correspond to better approximations.
Since the signal vector ${\bf f}$ in Fig.~1 is not ``smooth'', the 
inpainting approximation with the discrete Laplacian provides better 
results than the discrete biharmonic operator.
  

\def\arraystretch{0.6}
\setlength{\arraycolsep}{0.1em}
\def \imwidth {0.45\textwidth}

\begin{figure}[h]
\label{fig1}
\[
\begin{array}{c|cc}

& \mbox{Laplace operator} & \mbox{Biharmonic operator}\\[2mm]
\hline
\begin{sideways}\hspace{0.2cm}\parbox{4cm}{\centering
    \mbox{periodic boundary}\\ \mbox{conditions}\\[1mm]} \end{sideways}&
\includegraphics[width=\imwidth]
        {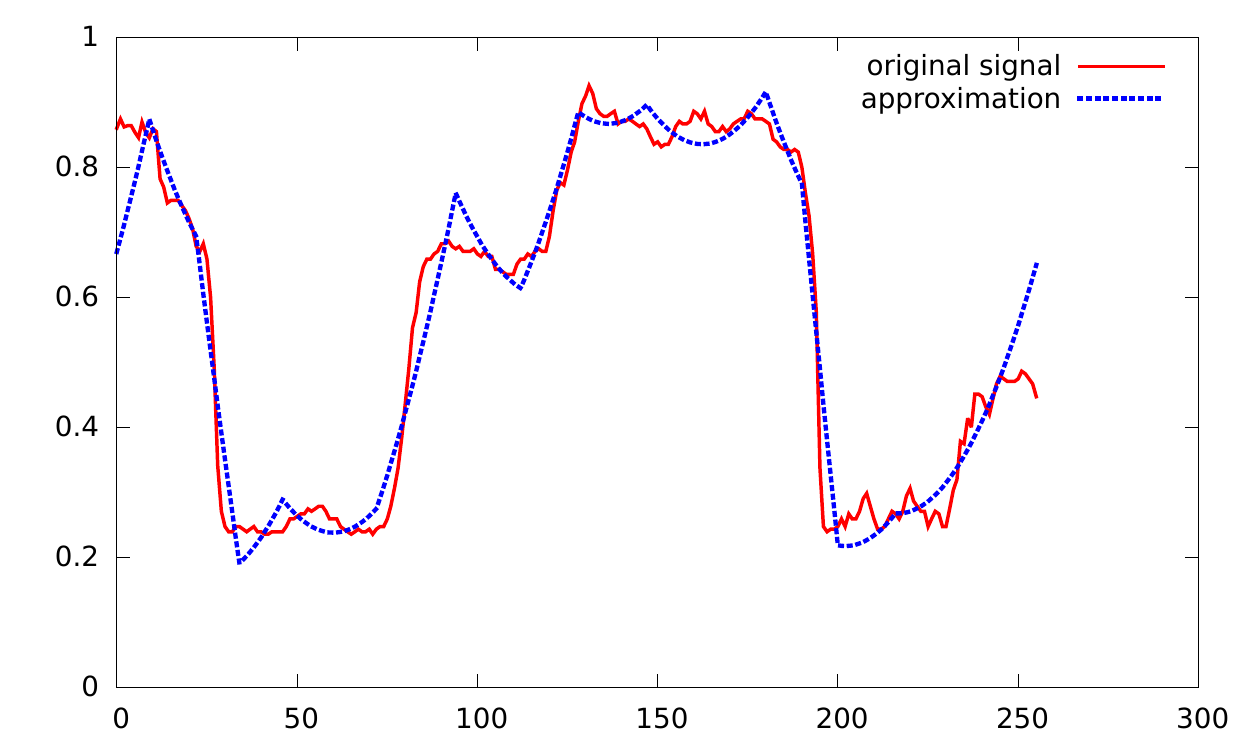} &
\includegraphics[width=\imwidth]
        {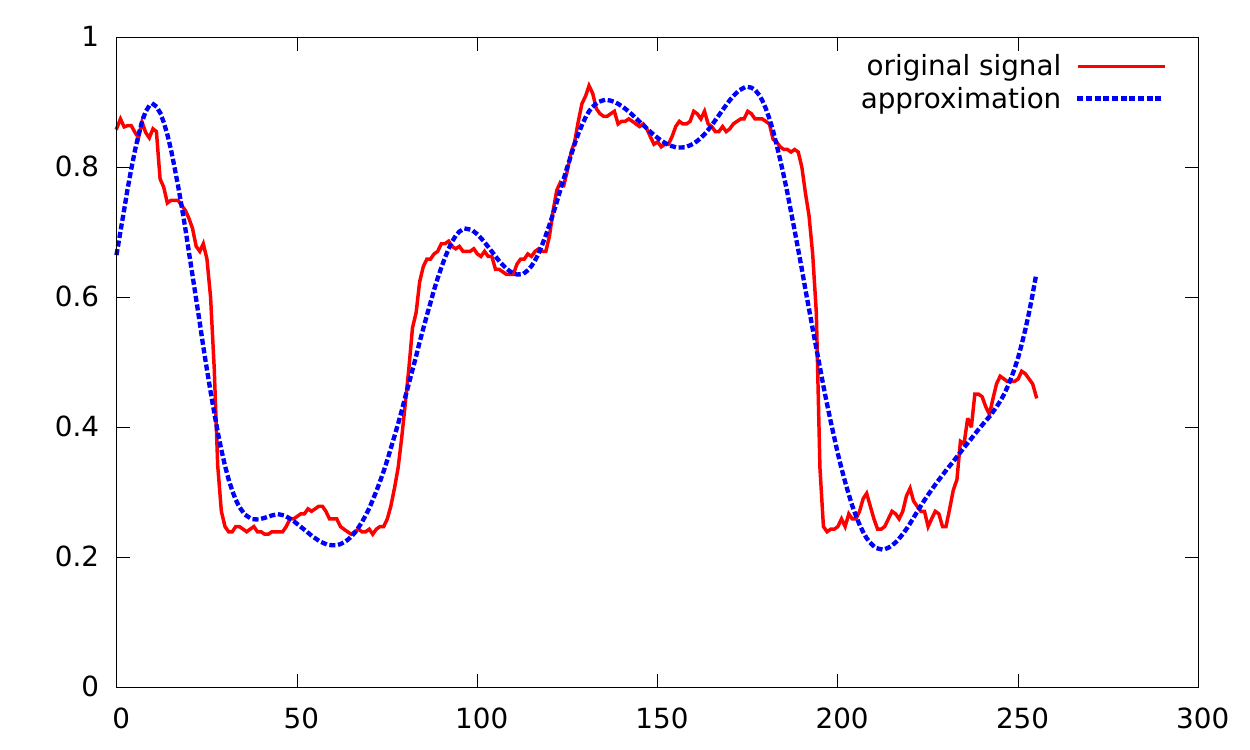}\\
& \mbox{PSNR: } 25.78 & \mbox{PSNR: } 24.39\\[2mm]

\begin{sideways}\hspace{0.05cm}\parbox{4.5cm}{\centering
        \mbox{homogeneous Neumann}\\ \mbox{boundary conditions}\\[1mm]} \end{sideways} &
\includegraphics[width=\imwidth]
        {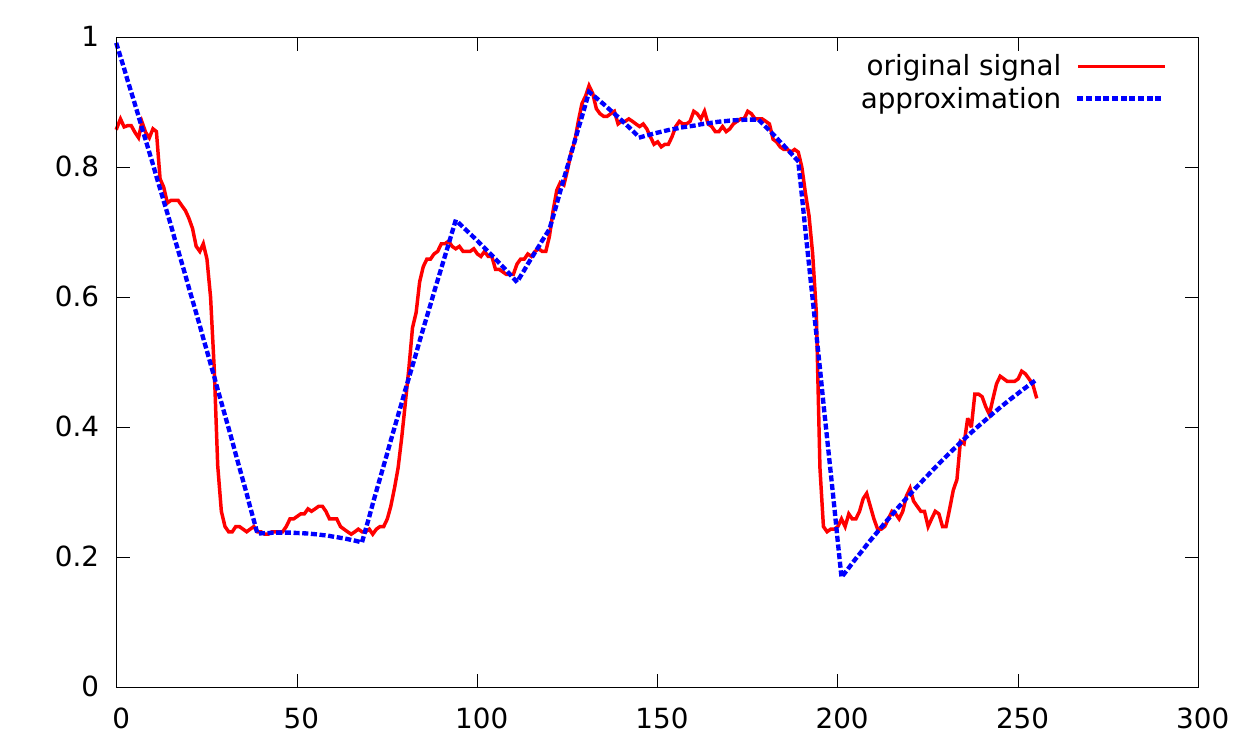} &
\includegraphics[width=\imwidth]
        {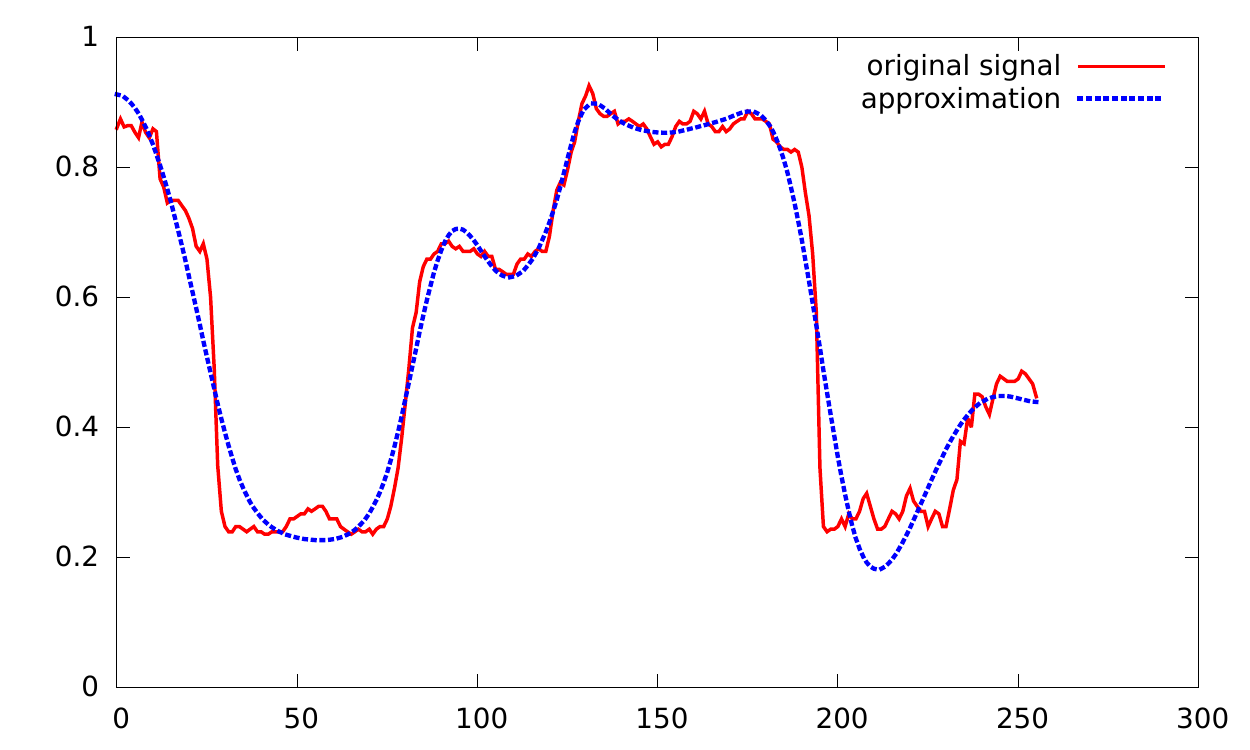} \\
& \mbox{PSNR: } 25.58 & \mbox{PSNR: } 25.30
\end{array}
\]
{\footnotesize {\bf Figure 1.} Sparse approximation of a vector ${\bf f} \in 
{\IR}^{256}$ (red line) by a vector {\bf u} (blue line) that is a linear 
combination of 13 columns of the pseudo-inverse ${\bf L}^+$ according to 
Algorithm 4.1. For ${\bf L}^+$, two different operators and two 
boundary conditions are considered.}
\end{figure}


\medskip
In the procedure described above, we have not implemented the 
``solvability condition''
\begin{equation} \label{sc}
 {\bf 1}^{T} {\bf g} =0 
\end{equation}
that we had derived from 
${\bf L} {\bf u} = {\bf g}$ and ${\bf 1}^{T} {\bf L} = {\bf 0}^{T}$ in 
Section 3. Therefore, the obtained  approximation results in Figure 1 
(top) for the Laplacian case appear to have piecewise quadratic form, 
since the columns of ${\bf A}^{+}$ and ${\bf B}^{+}$ in Theorems \ref{theo2} 
and \ref{theo3} are quadratic in $j$. In order to solve the original 
inpainting problem, we need also to incorporate the  constraint 
${\bf 1}^{T} {\bf g} =0$. This can be done by employing the additional 
equation $\sum_{i=1}^{|\Gamma|} g_{k_{i}} =0$ in the  last iteration 
step of the OMP algorithm when  updating the coefficients $g_{k_{1}}, 
\ldots , g_{k_{|\Gamma|}}$  by solving the least squares minimization problem.

 
\begin{figure}[h]
\label{fig2}
\[
\begin{array}{c|cc}

& \mbox{Laplace operator} & \mbox{Biharmonic operator}\\[2mm]
\hline
\begin{sideways}\hspace{0.2cm}\parbox{4cm}{\centering
    \mbox{periodic boundary}\\ \mbox{conditions}\\[1mm]} \end{sideways}&
\includegraphics[width=\imwidth]
        {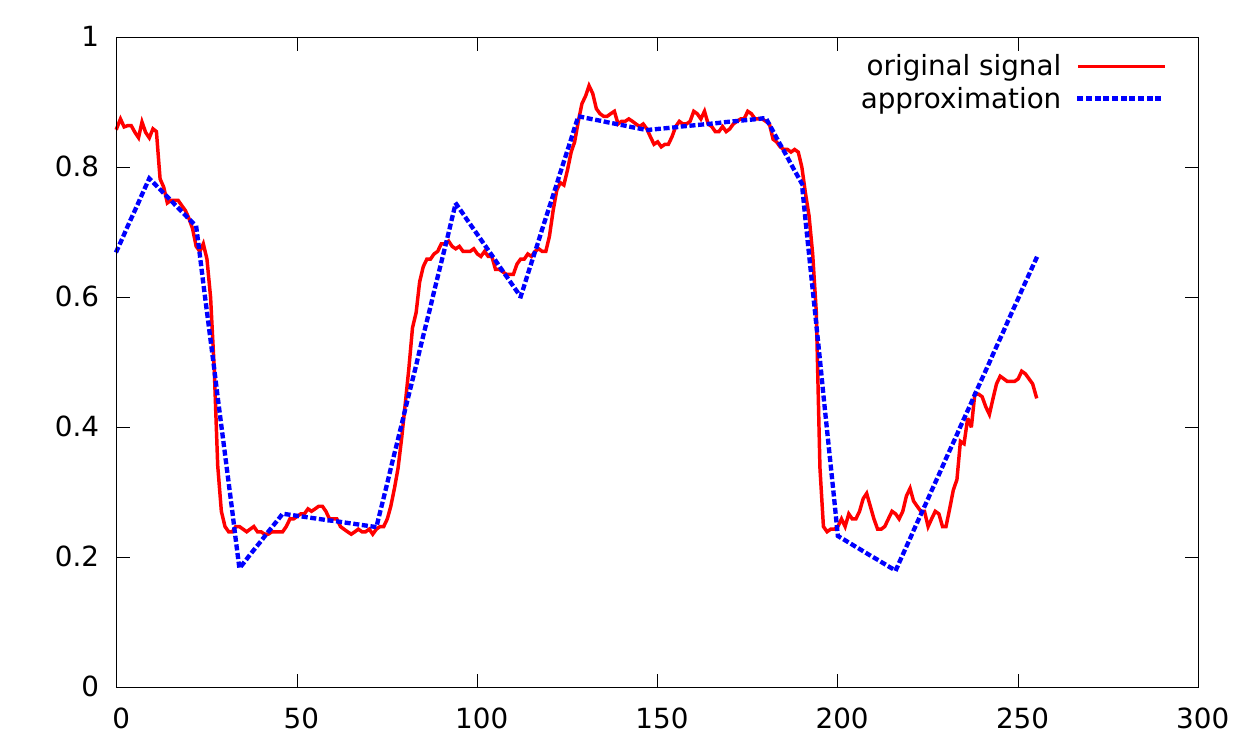} &
\includegraphics[width=\imwidth]
        {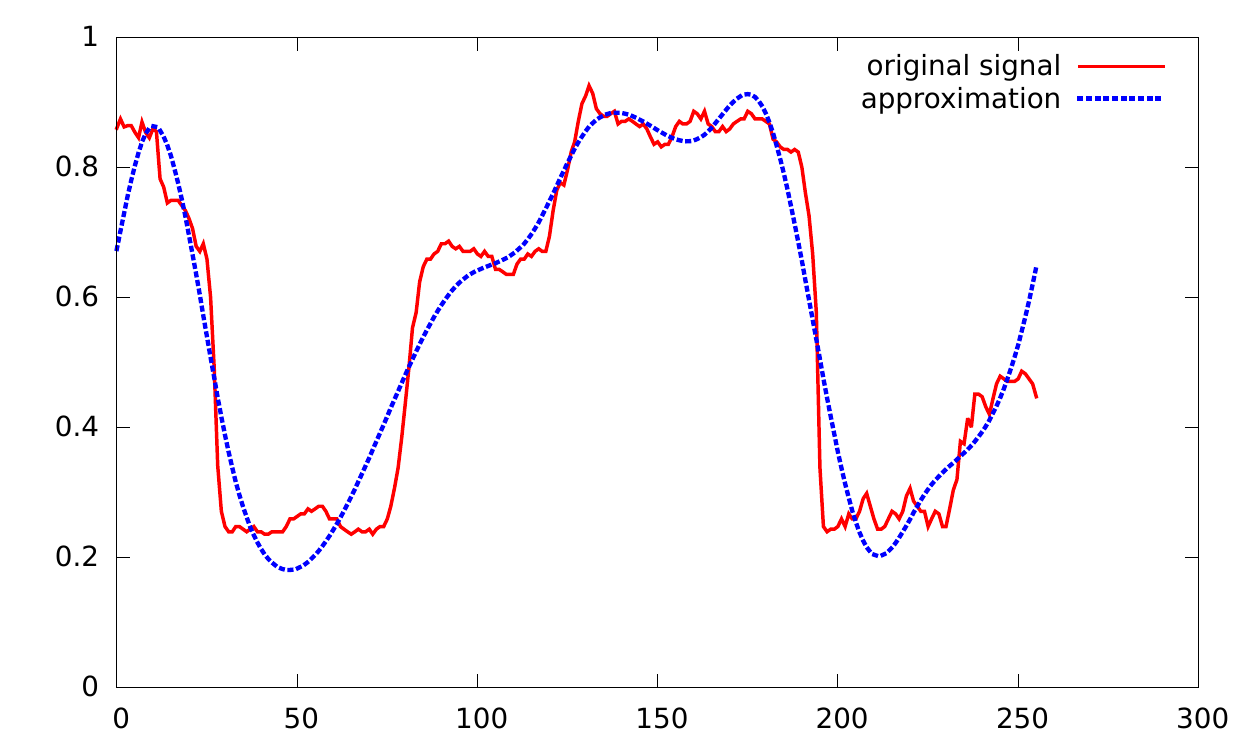}\\
& \mbox{PSNR: } 24.65 & \mbox{PSNR: } 23.21\\[2mm]

\begin{sideways}\hspace{0.05cm}\parbox{4.5cm}{\centering
        \mbox{homogeneous Neumann}\\ \mbox{boundary conditions}\\[1mm]} \end{sideways} &
\includegraphics[width=\imwidth]
        {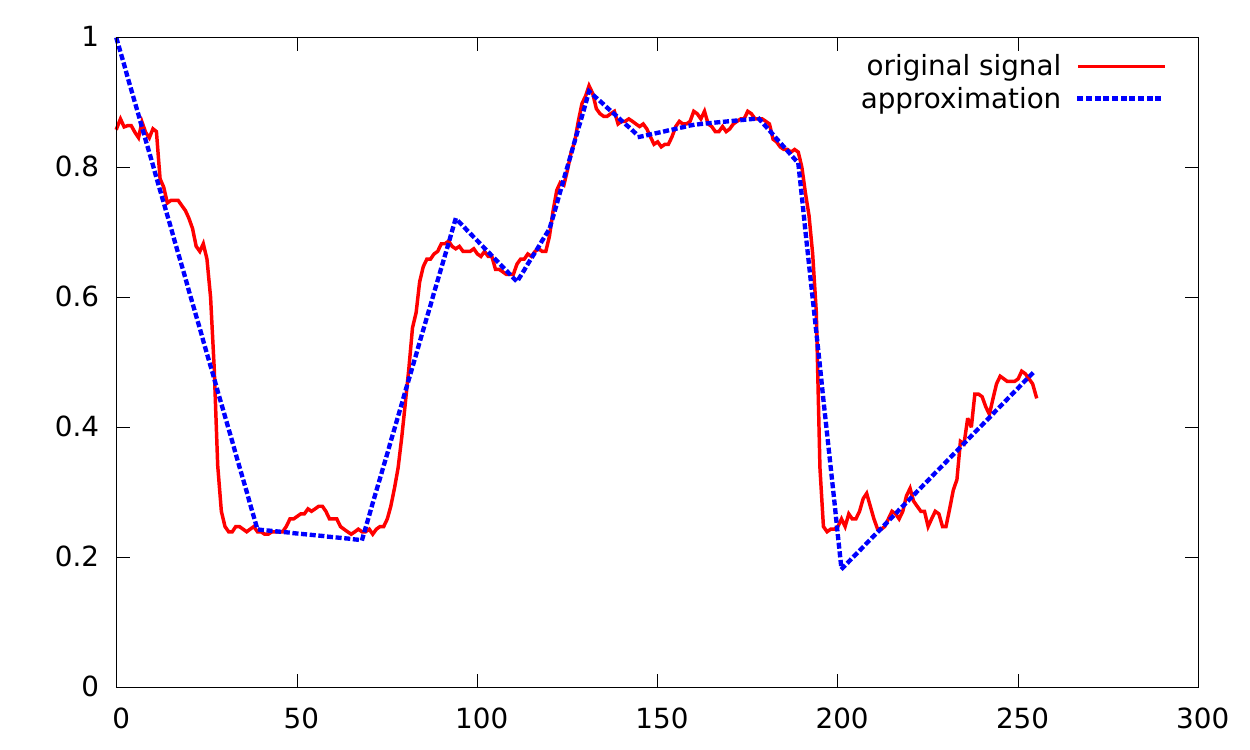} &
\includegraphics[width=\imwidth]
        {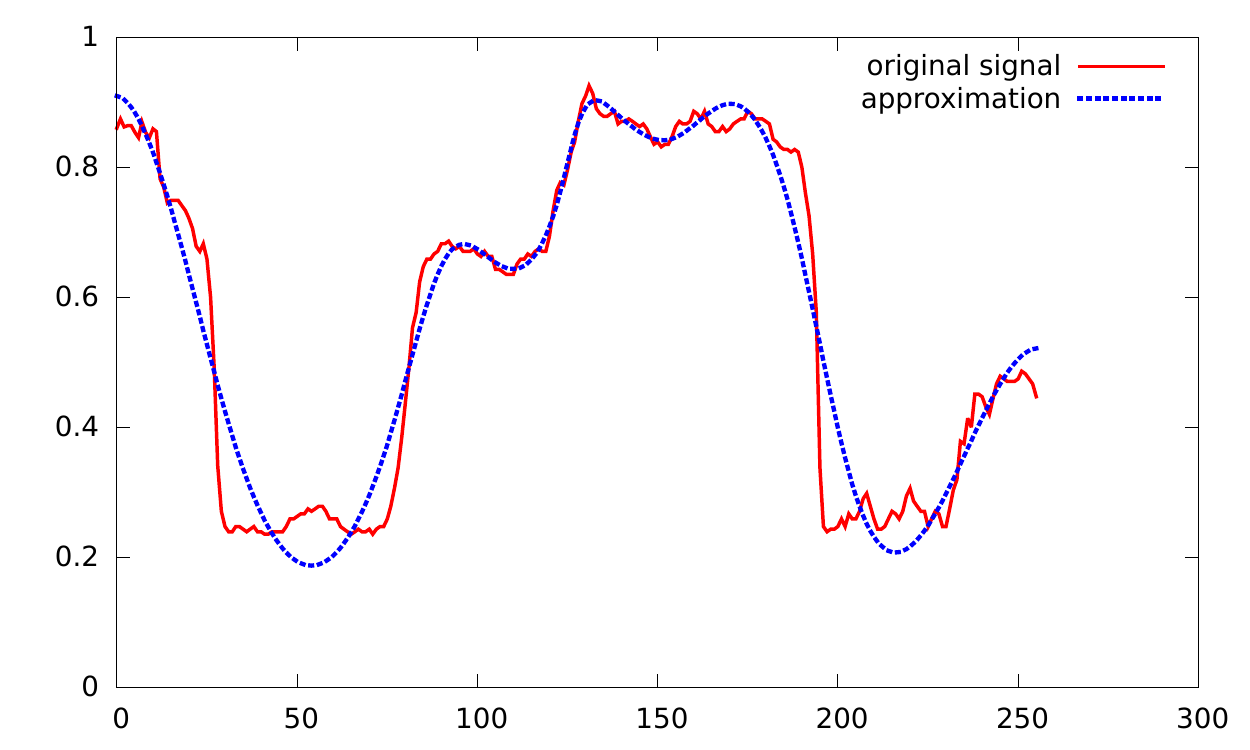} \\
& \mbox{PSNR: } 25.55 & \mbox{PSNR: } 24.16
\end{array}
\]
{\footnotesize {\bf Figure 2.} Sparse approximation of a vector ${\bf f} \in 
{\IR}^{256}$ (red line) by a vector {\bf u} (blue line) being a solution 
of the inpainting problem with different operators and boundary conditions. 
In contrast to Fig.~1, the solvability condition (\ref{sc}) has been taken 
into account.}
\end{figure}


\medskip
The results of the constrained sparse approximation, i.e., employing also 
the solvability condition,  for the four cases ${\bf L} \in \{ {\bf A}, 
\, {\bf B},  \, - {\bf A}^{2}, \, - {\bf B}^{2} \}$ are presented in 
Figure 2. Here we have used the same vector ${\bf f} \in {\IR}^{256}$ 
and $|\Gamma|=13$ as in the first experiment. In this case, we observe 
that the inpainting approach with the Laplacian provides 
a piecewise linear approximation.

\medskip
Both figures illustrate that the choice of the boundary conditions can have 
a fairly strong impact on the approximation result, even in locations that 
are far away from the boundaries. In 3 out of 4 cases, homogeneous Neumann 
boundary conditions gave better approximations than periodic ones. This
is to be expected for a non-periodic real-world signal. 


\end{document}